\documentclass{amsart}
\usepackage{amscd}
\usepackage{amssymb}
\usepackage{amsmath}
\usepackage{amsthm}
\usepackage{pifont}
\usepackage[cmtip,all]{xypic}

\theoremstyle{plain}
\newtheorem{Lem}{Lemma}[section]

\newtheorem{Cor}[Lem]{Corollary}
\newtheorem{Thm}[Lem]{Theorem}

{\theoremstyle{definition} 

\newtheorem{Rk}[Lem]{Remark}
\newtheorem{Def}[Lem]{Definition}}

{\theoremstyle{definition} 

}

\newcommand{\zig}{\addtocounter{Lem}{1}\tag{\theLem}}

\def\:{\colon}
\DeclareMathOperator*{\holim}{holim}
\DeclareMathOperator*{\1lim}{lim^1}
\DeclareMathOperator*{\lims}{lim^s}

\begin{document}

\title[The $E_2$-term of the descent spectral sequence]{The 
$E_2$-term of the descent spectral sequence for continuous 
$G$-spectra}

\author[Daniel G. Davis]{Daniel G. Davis$\sp 1$} 
\date{June 21, 2006}

\begin{abstract}
Given a profinite group $G$ with finite virtual cohomological dimension, 
let $\{X_i\}$ be a tower of discrete $G$-spectra, each of which is 
fibrant as a spectrum, so that $X=\holim_i X_i$ is a continuous $G$-spectrum, 
with homotopy fixed point spectrum $X^{hG}$. The $E_2$-term of the descent 
spectral sequence for $\pi_\ast(X^{hG})$ 
cannot always be expressed as continuous 
cohomology. However, we show that the $E_2$-term is always built out of 
a certain complex of spectra, that, in the context of abelian groups, is 
used to compute the continuous cochain cohomology of $G$ with coefficients 
in $\lim_i M_i$, where $\{M_i\}$ is a tower of discrete $G$-modules.
\end{abstract}

\maketitle
\footnotetext[1]{The author was supported by an NSF grant. Most of this 
paper was written during a visit to the Institut Mittag-Leffler (Djursholm, 
Sweden).}

\begin{section}{Introduction}
\par
In this note, $G$ always denotes a profinite group. Let $H_c^*(G;M)$ denote 
the 
continuous cohomology of $G$ with coefficients in the discrete $G$-module $M$. 
This cohomology is defined as the right derived functors of $G$-fixed points. 
Then we always assume that 
$G$ has finite virtual cohomological dimension; that is, there exists an open 
subgroup $H$ and a non-negative integer $m$, such that $H^s_c(H;M)=0$, for all 
discrete $H$-modules $M$ and all $s \geq m$.
\par
All of our spectra are Bousfield-Friedlander spectra of simplicial sets. 
In particular, a {\em discrete $G$-spectrum} is a $G$-spectrum such that each 
simplicial set $X_k$ is a simplicial object in the category of discrete 
$G$-sets (thus, for any $l \geq 0$, the action map on the 
$l$-simplices, $G \times (X_k)_l 
\rightarrow (X_k)_l$, is continuous when $(X_k)_l$ is regarded as a discrete 
space). The category of discrete $G$-spectra, with morphisms being 
$G$-equivariant maps of spectra, is denoted by $\mathrm{Spt}_G$. 
\par
Discrete $G$-spectra are considered in more detail in \cite{cts}, which shows 
(see \cite[Theorem 3.6]{cts}) 
that $\mathrm{Spt}_G$ is a model category, where a morphism 
$f$ in $\mathrm{Spt}_G$ is a weak equivalence (cofibration) if and only if 
$f$ is a weak equivalence (cofibration) in $\mathrm{Spt}$, the category of 
spectra. Given a discrete $G$-spectrum $X$, the {\em homotopy fixed point 
spectrum} $X^{hG}$ is obtained as the total right derived functor of fixed 
points: $X^{hG} = (X_{f,G})^G,$ where $X \rightarrow X_{f,G}$ is a trivial 
cofibration and $X_{f,G}$ is fibrant, all in 
$\mathrm{Spt}_G$.  
\par
Let $X_0 \leftarrow 
X_1 \leftarrow X_2 \leftarrow \cdots$ be a tower of discrete $G$-spectra, 
such that each $X_i$ is a fibrant spectrum. As explained in 
\cite[Lemma 4.4]{cts}, there exists a tower $\{X_i'\}$ of 
discrete $G$-spectra, such that there are weak equivalences 
\[\holim_i X_i \overset{\simeq}{\longrightarrow} \holim_i X_i' 
\overset{\simeq}{\longleftarrow} \lim_i X_i'.\] (In this paper, 
$\mathrm{holim}$ always denotes the version of the homotopy limit of spectra 
that is constructed levelwise in the category of simplicial sets, as defined 
in \cite{Bousfield/Kan} and \cite[5.6]{Thomason}.) Since the 
inverse limit of a tower of discrete $G$-sets is a topological $G$-space, 
and because $\holim_i X_i$ can be identified with $\lim_i X_i'$, in the eyes 
of homotopy, $\holim_i X_i$ is a 
{\em continuous $G$-spectrum}. Notice that, under this identification, 
the continuous $G$-action respects the 
topology of both $G$ and all the $X_i$ together. 
Continuous $G$-spectra and examples of such in 
chromatic stable homotopy theory are considered in \cite{cts}.
\par
Given the tower $\{X_i\}$ and the continuous $G$-spectrum $\holim_i X_i$, 
its {\em homotopy fixed point spectrum}, $(\holim_i X_i)^{hG}$, is 
defined to be $\holim_i (X_i)^{hG}.$ 
This construction is called homotopy fixed points 
because it is equivalent to the usual definition when $G$ is a finite group 
and it is the total right derived functor of fixed points in the appropriate 
sense (see \cite[Remark 8.4]{cts}).
\par
By \cite[Theorem 8.8]{cts}, thanks to the assumption of finite virtual 
cohomological dimension, there is a descent spectral sequence 
\begin{equation}\zig\label{dss}
E_2^{s,t} \Rightarrow \pi_{t-s}((\holim_i X_i)^{hG}),
\end{equation}
where 
\begin{equation}\zig\label{e2}
E_2^{s,t} = \pi^s\pi_t(\holim_i (\Gamma^\bullet_G(X_i)_{f,G})^G)
\end{equation}
(see the beginning of Section 2 for the meaning of $\Gamma_G^\bullet$), 
and, 
if the tower of abelian groups $\{\pi_t(X_i)\}$ satisfies the Mittag-Leffler 
condition for every integer $t$, then $E_2^{s,t} \cong 
H^s_\mathrm{cont}(G;\{\pi_t(X_i)\}),$ which is continuous cohomology in the 
sense of Jannsen. (This cohomology is obtained by taking the 
right derived functors of $\lim_i (-)^G$, a functor from towers of 
discrete $G$-modules to abelian groups; see \cite{Jannsen}.) 
\par
In expression (\ref{e2}), since $\pi_t(\holim_i (-))$ is not necessarily 
$\lim_i \pi_t(-)$, the $E_2$-term of descent spectral sequence (\ref{dss}), 
in general, can not be 
expressed as continuous cohomology, and, in general, it has no 
compact algebraic 
description. For example, as pointed out in \cite[Remark 8.10]{cts}, when 
$G=\{e\},$ \[H^s_\mathrm{cont}(\{e\};\{\pi_t(X_i)\}) = \lim_i \pi_t(X_i)\] and 
$E_2^{0,t} = \pi_t(\holim_i X_i)$, and these need not be isomorphic, due to 
the familiar $\1lim_i \pi_{t+1}(X_i)$ obstruction. 
However, in this note, we show that the $E_2$-term (\ref{e2}) 
can always 
be described in an interesting way. 
\par
In more detail, Theorem \ref{main1} gives a particular cochain complex 
$\mathcal{C}^\ast$ for computing 
the continuous cochain cohomology of $G$ 
for a topological $G$-module $\lim_i M_i,$ where $\{M_i\}$ is a tower of 
discrete $G$-modules. In Corollary \ref{main3}, we show 
that the $E_2$-term of (\ref{e2}) can always 
be given by taking the cohomology of 
the homotopy groups of the complex $\mathcal{C}^\ast$, where 
$\lim_i M_i$ is replaced 
by the continuous $G$-spectrum $\holim_i X_i$, in an appropriate sense. This 
presentation of the $E_2$-term shows that $E_2^{\ast,\ast}$ always takes into 
account the topology of the continuous $G$-spectrum, even when it cannot be
expressed as continuous cohomology.
\par
Section 3 of this note consists of a discussion of the construction of descent 
spectral sequence (\ref{dss}). Section 5 explains why two other potentially 
plausible interpretations of (\ref{e2}) fail to work.
\vspace{.1in}
\par
\noindent
\textbf{Acknowledgements.} I am grateful to the referee for helpful comments 
and for pointing out that the statement of the main result could be 
simplified. Also, I thank Paul Goerss and Halvard Fausk for useful remarks.
\end{section}
\begin{section}{The pro-discrete cochain complex and continuous cohomology}
\par
We begin this section with some terminology. 
If $\mathcal{C}$ is a category, then $\mathbf{tow}(\mathcal{C})$ is the 
category of towers \[C_0 \leftarrow C_1 \leftarrow C_2 \leftarrow \cdots\] in 
$\mathcal{C}$. The morphisms $\{f_i\}$ are natural transformations such that 
each $f_i$ is a morphism in $\mathcal{C}$. 
In this note, we will be working with 
$\mathbf{tow}(\mathrm{DMod}(G))$, where $\mathrm{DMod}(G)$ is the category of 
discrete $G$-modules, and $\mathbf{tow}(\mathrm{Spt}_G)$.
\par
If $A$ is an abelian group with the discrete topology, let 
$\mathrm{Map}_c(G,A)$ be the abelian group of continuous maps from $G$ to $A$. 
If $X$ is a spectrum, one can also define $\mathrm{Map}_c(G,X),$ where 
the $l$-simplices of the $k$th simplicial set $(\mathrm{Map}_c(G,X)_k)_l$ are 
given by $\mathrm{Map}_c(G,(X_k)_l),$ where $(X_k)_l$ is given the discrete 
topology. 
\par
Consider the functor 
\[\Gamma_G \colon \mathrm{Spt}_G \rightarrow \mathrm{Spt}_G, \ \ \ 
X \mapsto \Gamma_G(X) = \mathrm{Map}_c(G,X),\] where the action of $G$ on 
$\mathrm{Map}_c(G,X)$ 
is induced on the level of sets 
by $(g \cdot f)(g') = f(g'g),$ for $g, g' \in G$ and 
$f \in \mathrm{Map}_c(G,(X_k)_l),$ for each $k, l \geq 0.$ 
As explained 
in \cite[Definition 7.1]{cts}, the functor $\Gamma_G$ 
forms a triple and there is 
a cosimplicial discrete $G$-spectrum $\Gamma^\bullet_G X.$ Also, it is clear 
that $\Gamma_G \: \mathrm{DMod}(G) \rightarrow \mathrm{DMod}(G)$ can be 
defined as above, so that, given a discrete $G$-module $M$, 
$\Gamma^\bullet_G M$ is a cosimplicial discrete $G$-module.
\par
We do not claim any originality in the following definition.
\begin{Def}
Let $\{X_i\}$ be an object in $\mathbf{tow}(\mathrm{DMod}(G))$ or 
in $\mathbf{tow}(\mathrm{Spt}_G)$. 
Then the {\em pro-discrete cochain complex} is defined to be the complex 
\[\mathcal{C}^\ast(G;\{X_i\}) = 
\lim_i (\Gamma^\ast_G X_i)^G,\] where 
$\Gamma^\ast_G X_i$ is the canonical complex associated to 
$\Gamma^\bullet_G X_i,$ where, if $\{X_i\}$ is in 
$\mathbf{tow}(\mathrm{Spt}_G)$, then the complex lives in the stable 
homotopy category. The pro-discrete cochain complex is a complex of 
abelian groups or spectra, respectively, and the limit and colimit are both 
formed in abelian groups or in spectra (not the stable homotopy 
category), respectively.
\end{Def}
\par
Let $M$ be any topological $G$-module (that is, an abelian 
Hausdorff topological group that is a $G$-module, with a continuous 
$G$-action). Then the {\em 
continuous cochain cohomology 
of $G$ with coefficients in $M$}, $H^*_\mathrm{cts}(G;M)$, 
is the cohomology of a cochain complex that has the form 
\begin{equation}\zig\label{cplx}
M \rightarrow \mathrm{Map}_c(G,M) \rightarrow \mathrm{Map}_c(G^2,M) 
\rightarrow \cdots
\end{equation} (see \cite[pg. ~106]{Neuk} for details). We note that, 
by \cite[Theorem (2.2)]{Jannsen}, if $\{M_i\}$ is a tower of discrete 
$G$-modules that satisfies the Mittag-Leffler condition, then 
\[H^s_\mathrm{cts}(G;\lim_i M_i) \cong H^s_\mathrm{cont}(G;\{M_i\}),\] 
for all $s \geq 0$, but, 
in general, these two versions of continuous cohomology need not be 
isomorphic. Also, if $M$ is a discrete $G$-module, then 
$H^s_\mathrm{cts}(G;M) = H^s_c(G;M).$
\par
Now we show that the pro-discrete cochain complex can be used to compute 
continuous cochain cohomology.
\begin{Thm}\label{main1}
If $\{M_i\}$ is a tower of discrete $G$-modules, then 
\[H^s_\mathrm{cts}(G;\lim_i M_i) \cong H^s[\mathcal{C}^\ast(G;\{M_i\})].\]
\end{Thm} 
\begin{proof}
As explained in \cite[pg. ~106]{Neuk}, for a topological $G$-module $M$, 
the chain complex in (\ref{cplx}) is 
defined by taking 
the $G$-fixed points of the complex 
\begin{equation}\zig\label{cplx2}
X^\ast(G;M) = [\mathrm{Map}_c(G,M) \rightarrow \mathrm{Map}_c(G^2,M) 
\rightarrow \cdots], 
\end{equation}
where $X^n(G;M) = \mathrm{Map}_c(G^{n+1},M)$ has a 
$G$-action that is defined by 
\[(g\cdot f)(g_1,...,g_{n+1}) = g \cdot f(g^{-1}g_1, ..., g^{-1}g_{n+1}).\]
\par
Now let $M$ be a discrete $G$-module. Then 
it is a standard fact that the cochain complex $(X^\ast(G,M))^G$ is 
naturally isomorphic as a 
complex to the cochain 
complex $(\Gamma^\ast_G M)^G.$ This isomorphism uses the fact that the 
abelian group of $n$-cochains 
of $(\Gamma^\ast_G M)^G$ is isomorphic to $\mathrm{Map}_c(G^{n+1},M)^G,$ 
where $\mathrm{Map}_c(G^{n+1},M)$ has a $G$-action that 
is given by \[(g\cdot f)(g_1,g_2,g_3,...,g_{n+1}) 
= f(g_1g,g_2,g_3,...,g_{n+1}).\]
\par
Since
\[(X^n(G;\lim_i M_i))^G \cong 
\lim_i ((X^n(G;M_i))^G) \cong \lim_i (\Gamma^{n+1}_G M_i)^G,\] 
we have:
\[H^s_\mathrm{cts}(G;\lim_i M_i) = H^s[(X^\ast(G;\lim_i M_i))^G] 
= H^s[\lim_i (\Gamma^\ast_G M_i)^G],\] where we used the aforementioned fact 
that $(X^\ast(G,M_i))^G$ and $(\Gamma^\ast_G M_i)^G$ are naturally isomorphic 
cochain complexes.
\end{proof}
\end{section}
\begin{section}{Constructing the descent spectral sequence}
\par
In this section, we review how descent spectral sequence (\ref{dss}) is 
constructed and we compare it with a spectral sequence whose 
$E_2$-term is always Jannsen's continuous cohomology.
\par
Given a tower $\{X_i\}$ of discrete $G$-spectra, such that each $X_i$ is 
a fibrant spectrum, by \cite[Remark 7.8, Definition 8.1]{cts}, 
\[(\holim_i X_i)^{hG} = 
\holim_i \holim_\Delta (\Gamma^\bullet_G (X_i)_{f,G})^G.\] Thus, 
\[(\holim_i X_i)^{hG} \cong 
\holim_\Delta (\holim_i (\Gamma^\bullet_G (X_i)_{f,G})^G),\] and 
descent spectral sequence (\ref{dss}) is the conditionally 
convergent homotopy spectral 
sequence 
\[\lims_\Delta \pi_t(\holim_i (\Gamma^\bullet_G (X_i)_{f,G})^G) 
\Rightarrow 
\pi_{t-s}(\holim_\Delta (\holim_i (\Gamma^\bullet_G (X_i)_{f,G})^G))\] 
(see \cite[Theorem 8.8]{cts}).
\par
In the above context, there is another spectral sequence that is natural 
to consider. Since 
\[(\holim_i X_i)^{hG} \cong
\holim_{\Delta \times \{i\}} (\Gamma^\bullet_G (X_i)_{f,G})^G,\] there is 
a conditionally convergent homotopy spectral sequence 
\[\lims_{\Delta \times \{i\}} \pi_t((\Gamma^\bullet_G(X_i)_{f,G})^G) 
\Rightarrow 
\pi_{t-s}(\holim_{\Delta \times \{i\}} (\Gamma^\bullet_G (X_i)_{f,G})^G) 
\cong \pi_{t-s}((\holim_i X_i)^{hG}),\] 
such that 
\[\lims_{\Delta \times \{i\}} \pi_t((\Gamma^\bullet_G(X_i)_{f,G})^G) 
\cong H^s_\mathrm{cont}(G;\{\pi_t(X_i)\})\] 
(see \cite[Proposition 3.1.2]{Geisser}). This spectral sequence is 
closely related to the $\ell$-adic descent spectral sequence of algebraic 
$K$-theory (see \cite{Thomason}, \cite{Mitchell}).
\par
We see that we have two spectral sequences with abutment 
$\pi_\ast((\holim_i X_i)^{hG})$. As pointed out in the Introduction, 
if the tower $\{\pi_t(X_i)\}$ satisfies the Mittag-Leffler condition for 
every integer $t$, then the spectral sequences have isomorphic $E_2$-terms. 
However, as stated in the Introduction, 
the case $G=\{e\}$ shows that these two spectral sequences can have 
different $E_2$-terms, so that the spectral sequences can be different 
from each other. 
\par
Since the second spectral sequence described above has an $E_2$-term that 
always has a nice algebraic description, it is natural to ask what is the 
value of descent spectral sequence (\ref{dss}). We will see that, because 
the descent spectral sequence is the homotopy spectral sequence of a 
cosimplicial spectrum, in certain cases it can be compared with an 
Adams-type spectral sequence that is strongly convergent. 
\par
Let $n \geq 1$ and let $p$ be a prime. 
Let $K(n)$ be the $n$th Morava $K$-theory spectrum with $K(n)_\ast = 
\mathbb{F}_p[v_n^{\pm 1}],$ where the degree of $v_n$ is $2(p^n-1)$. Also, 
let $E_n$ denote the Lubin-Tate spectrum, where $E_{n \ast} = 
W(\mathbb{F}_{p^n})[[u_1, ..., u_{n-1}]][u^{\pm 1}],$ where the degree of 
$u$ is $-2$, and the complete power series ring over the Witt vectors is in 
degree zero. 
\par
Let $Z$ be a $K(n)$-local spectrum and suppose that there is an 
augmentation $Z \rightarrow \holim_i (\Gamma^\bullet_G (X_i)_{f,G})^G$. If 
the associated complex of spectra
\[\ast \rightarrow Z \rightarrow \holim_i \mathrm{Map}_c(G,(X_i)_{f,G})^G 
\rightarrow \holim_i \mathrm{Map}_c(G,\mathrm{Map}_c(G,(X_i)_{f,G}))^G 
\rightarrow \cdots \] is a $K(n)$-local $E_n$-resolution of $Z$ (for 
the definition of this, see \cite[Appendix A]{DH}), then descent spectral 
sequence (\ref{dss}) is isomorphic to the strongly convergent 
$K(n)$-local $E_n$-Adams 
spectral sequence with abutment $\pi_\ast(Z)$ (see \cite[Proposition A.5, 
Corollary A.8]{DH}). Thus, the descent spectral sequence is strongly 
convergent and the map $Z \rightarrow \holim_\Delta 
\holim_i (\Gamma^\bullet_G (X_i)_{f,G})^G \cong (\holim_i X_i)^{hG}$ is a weak 
equivalence (\cite[Corollary A.8]{DH}). 
\par
In this way, in 
\cite[Chapter 10]{thesis}, the author showed that, given a finite spectrum 
$X$, the descent spectral sequence for $\pi_\ast((E_n \wedge X)^{hG})$ 
is strongly convergent and isomorphic to the $K(n)$-local $E_n$-Adams 
spectral sequence with abutment $\pi_\ast(E_n^{dhG} \wedge X)$, where $G$ 
is a closed subgroup of the extended Morava stabilizer group $G_n$ and 
$E_n^{dhG}$ is the spectrum constructed by Devinatz and Hopkins in 
\cite{DH} ($E_n^{dhG}$ is denoted by $E_n^{hG}$ in \cite{DH}).
\end{section}
\begin{section}{The $E_2$-term and the pro-discrete cochain 
complex}
\par
In this section, we show that the $E_2$-term of (\ref{e2}) can be built out 
of the same complex that computes continuous cochain cohomology. More 
precisely, given a continuous $G$-spectrum $\holim_i X_i$, there exists a 
tower $\{X_i'\}$ of discrete $G$-spectra, such that 
\begin{equation}\zig\label{neato}
E_2^{s,t} \cong H^s[\pi_t(\mathcal{C}^\ast(G;\{X_i'\}))].\end{equation} 
\par
We find the expression on the right-hand side in (\ref{neato}) interesting 
for the following reason. 
The homotopy fixed point spectrum is defined with respect to a 
continuous action of $G$ on the spectrum. Thus, homotopy fixed points take 
into account the topology of the spectrum. Similarly, since the $E_2$-term 
is built out of the pro-discrete cochain complex of spectra, the 
$E_2$-term is always taking into 
account the topology of the spectrum.
\par
By \cite[VI, Proposition 1.3]{GJ}, 
$\mathbf{tow}(\mathrm{Spt}_G)$ is a model category, where $\{f_i\}$ is a 
weak equivalence (cofibration) if and only if each $f_i$ is a weak equivalence 
(cofibration) in $\mathrm{Spt}_G$. 
\begin{Thm}\label{main2}
The $E_2$-term (\ref{e2}) of descent spectral sequence (\ref{dss}) 
has the form
\begin{equation}\zig\label{nice}
E_2^{s,t} \cong \pi^s\pi_t(\lim_i (\Gamma^\bullet_G X_i')^G),
\end{equation}
where $\{X_i\} \rightarrow
\{X_i'\}$ is a trivial cofibration with $\{X_i'\}$ fibrant, all in
$\mathbf{tow}(\mathrm{Spt}_G)$.
\end{Thm}
\begin{proof}
Let $\{X_i'\}$ be as stated in the theorem. 
By \cite[VI, Remark 1.5]{GJ}, each $X_i'$ is
fibrant and each map $X_i' \rightarrow X_{i-1}'$ is a fibration, all
in $\mathrm{Spt}_G$.
\par
For any $k \geq 0$, we consider the expression
\[\holim_i ((\Gamma^\bullet_G X_i')^G)^k =
\holim_i (\mathrm{Map}_c(G,\mathrm{Map}_c(G, \cdots
,\mathrm{Map}_c(G, X_i') \cdots )))^G,\] where $\mathrm{Map}_c(G,-)$
appears $k+1$ times. By \cite[Section 3]{cts}, the forgetful functor 
$U \: \mathrm{Spt}_G
\rightarrow \mathrm{Spt},$ 
$\mathrm{Map}_c(G,-) \: \mathrm{Spt} \rightarrow \mathrm{Spt}_G,$ where 
$\mathrm{Map}_c(G,X) = \Gamma_G(X)$, and
the functor $(-)^G \: \mathrm{Spt}_G \rightarrow 
\mathrm{Spt}$ all preserve
fibrations. Thus, $\{X_i'\}$ is a tower of fibrations of fibrant spectra,
all in $\mathrm{Spt}.$ 
This implies that $\{\mathrm{Map}_c(G,X_i')\}$ is a tower of
fibrations of fibrant spectra, in $\mathrm{Spt}_G$, and hence, in 
$\mathrm{Spt}.$ By
iteration, \[\{\mathrm{Map}_c(G,\mathrm{Map}_c(G, \cdots,
\mathrm{Map}_c(G, X_i') \cdots ))\}\] is a tower of fibrations of
fibrant spectra, in $\mathrm{Spt}_G$, so that
\[\{(\mathrm{Map}_c(G,\mathrm{Map}_c(G, \cdots,
\mathrm{Map}_c(G, X_i') \cdots )))^G\}\] is a tower of fibrations of
fibrant spectra in $\mathrm{Spt}.$ Therefore, the canonical map \[\lim_i
((\Gamma^\bullet_G X_i')^G)^k \rightarrow \holim_i
((\Gamma^\bullet_G X_i')^G)^k\] is a weak equivalence.
\par
Since $\{((\Gamma^\bullet_G X_i)^G)^k\}$ and
$\{((\Gamma^\bullet_G (X_i)_{f,G})^G)^k\}$ are towers of fibrant spectra,
there is a zigzag of weak equivalences
\[\lim_i ((\Gamma^\bullet X_i')^G)^k \rightarrow \holim_i
((\Gamma^\bullet X_i')^G)^k \leftarrow \holim_i
((\Gamma^\bullet X_i)^G)^k \rightarrow \holim_i
((\Gamma^\bullet (X_i)_{f,G})^G)^k,\] where $\Gamma = \Gamma_G$. 
This zigzag of weak
equivalences implies that
\[\pi^s\pi_t(\lim_i
((\Gamma^\bullet_G X_i')^G)) \cong \pi^s\pi_t(\holim_i
((\Gamma^\bullet_G (X_i)_{f,G})^G)).\]
\end{proof}
\begin{Cor}\label{main3}
Let $\{X_i'\}$ be as in Theorem \ref{main2}. Then there is an isomorphism 
\[E_2^{s,t} \cong H^s[\pi_t(\mathcal{C}^\ast(G;\{X_i'\}))],\] where 
$E_2^{s,t}$ is the $E_2$-term of (\ref{e2}).
\end{Cor}
\begin{Rk}
By Theorem \ref{main1}, $H^s[\mathcal{C}^\ast(G;\{\pi_t(X_i')\})] 
\cong H^s_\mathrm{cts}(G;\lim_i \pi_t(X_i)).$
\end{Rk}
\end{section}
\begin{section}{The failure of other possible descriptions of the $E_2$-term}
After studying the expression in (\ref{nice}) further, one recalls that 
$\lim_i(-)^G$ is the functor used to define $H^s_\mathrm{cont}(G;-)$, and, if 
$M$ is any discrete $G$-module, then 
\[0 \rightarrow M \rightarrow \Gamma^\ast_G M\] 
is a $(-)^G$-acyclic resolution of $M$, so that 
$H^s[(\Gamma^\ast_G M)^G] = H^s_c(G;M).$
\par
Let $\{M_i\}$ be a tower of discrete $G$-modules. If 
\[\{0\} \rightarrow \{M_i\} \rightarrow \{\Gamma^\ast_G M_i\}\] is a 
$\lim_i(-)^G$-acyclic resolution of $\{M_i\}$ in 
$\mathbf{tow}(\mathrm{DMod}(G))$, then 
\[H^s[(\lim_i(-)^G)(\{\Gamma^\ast_G M_i\})] = H^s_\mathrm{cont}(G;\{M_i\}).\]
This would imply that $E_2^{s,t} \cong 
H^s[\pi_t(\lim_i (\Gamma^\ast_G X_i')^G)]$ is computed by taking the 
cohomology of the homotopy groups of a complex of spectra in the 
stable homotopy category, that, in the 
context of abelian groups, computes continuous cohomology. This would be 
an interesting presentation of the $E_2$-term. 
\par
However, it is not hard to show that 
\[\{0\} \rightarrow \{M_i\} \rightarrow \{\Gamma^\ast_G M_i\}\] need not be a 
$\lim_i(-)^G$-acyclic resolution of $\{M_i\}$ in 
$\mathbf{tow}(\mathrm{DMod}(G))$, so that the above interpretation of the 
$E_2$-term does not work out. For example, 
by \cite[(2.1)]{Jannsen}, there is a short 
exact sequence 
\[0 \rightarrow \1lim_i H^{s-1}_c(G;\Gamma_G M_i) \rightarrow 
H^s_\mathrm{cont}(G;\{\Gamma_G M_i\}) 
\rightarrow \lim_i H^s_c(G;\Gamma_G M_i) \rightarrow 0,\] for each $s \geq 0$, 
where $H^{-1}_c(G;-)=0$. 
Therefore, when $s \geq 1$, $H^s_c(G;\Gamma_G M_i)=0$, so that, for all $s 
\geq 2$, 
$H^s_\mathrm{cont}(G;\{\Gamma_G M_i\}) = 0.$ 
But, the short exact sequence also implies 
that \[H^1_\mathrm{cont}(G;\{\Gamma_G M_i\}) 
\cong \1lim_i M_i,\] which need not 
vanish. Thus, $\{\Gamma_G M_i\}$, the first object in the complex 
$\{\Gamma^\ast_G M_i \}$, need not be $\lim_i(-)^G$-acyclic 
in $\mathbf{tow}(\mathrm{DMod}(G))$.
\par
Upon further consideration of the expression in (\ref{nice}), one notices 
that, for any $k,l,m \geq 0$, 
\[((\lim_i(\Gamma^{m+1}_G X_i')^G)_k)_l = 
\lim_i(\Gamma^{m+1}_G((X_i')_k)_l)^G \cong 
\mathrm{Map}_c(G^m, \lim_i((X_i')_k)_l)\]
is an isomorphism of sets. If one could promote this isomorphism to 
\begin{equation}\zig\label{iso}
\lim_i(\Gamma^{m+1}_G X_i')^G \cong 
\mathrm{Map}_c(G^m, \lim_i X_i'),
\end{equation} then one could use this to interpret 
the expression in (\ref{nice}) as being the cohomology of homotopy groups 
applied to the complex of continuous cochains with target (``coefficients'') 
the continuous 
$G$-spectrum $\lim_i X_i'$.
\par
But notice that, in this interpretation, 
the expression $\mathrm{Map}_c(G^m, \lim_i X_i')$ does not have 
the desired meaning. For isomorphism (\ref{iso}) to hold, $\lim_i X_i'$ must 
be a 
spectrum whose simplicial sets have simplices with the pro-discrete topology. 
But, as a Bousfield-Friedlander spectrum, in the construction 
$\mathrm{Map}_c(G^m, \lim_i X_i')$, $\lim_i X_i'$ consists of simplicial 
sets whose simplices all have the discrete topology, by default. 
This conflict means that this interpretation also fails to 
work. 

\end{section}
\bibliographystyle{amsplain}

\bibliography{biblio}

\end{document}